\documentclass{lms}
\usepackage{amscd,amssymb,amsfonts,latexsym,subfigure}
\usepackage[arrow,matrix,graph]{xy}

\newtheorem{theorem}{Theorem}[section]
\newtheorem{corollary}[theorem]{Corollary}
\newtheorem{lemma}[theorem]{Lemma}
\newtheorem{prop}[theorem]{Proposition}

\newnumbered{definition}[theorem]{Definition}
\newnumbered{example}[theorem]{Example}
\newnumbered{remark}[theorem]{Remark}
\newunnumbered{ack}{Acknowledgments}
\newunnumbered{claim}{Claim}

\newcommand{\N}{\mathbb{N}}
\newcommand{\Z}{\mathbb{Z}}
\newcommand{\Q}{\mathbb{Q}}

\newcommand{\C}{\mathbb{C}}

\newcommand{\TT}{\mathbb{T}}

\renewcommand{\k}{\Bbbk}
\newcommand{\HH}{{\mathfrak H}}
\newcommand{\B}{{\mathfrak B}}
\newcommand{\RR}{{\mathcal R}_1}
\newcommand{\VV}{{\mathcal V}_1}

\newcommand{\G}{\Gamma}
\newcommand{\E}{\mathsf{E}}
\newcommand{\V}{\mathsf{V}}
\newcommand{\T}{\mathsf{T}}
\newcommand{\W}{\mathsf{W}}

\providecommand{\rank}{\mathop{\rm rank}}
\providecommand{\gr}{\mathop{\rm gr}}
\providecommand{\im}{\mathop{\rm im}}
\providecommand{\coker}{\mathop{\rm coker}}
\providecommand{\codim}{\mathop{\rm codim}}
\providecommand{\id}{\mathop{\rm id}}
\providecommand{\ab}{\mathop{\rm ab}}
\providecommand{\Sym}{\mathop{\rm Sym}}

\providecommand{\Lie}{\mathop{\rm Lie}}

\providecommand{\Hom}{\mathop{\rm Hom}}

\providecommand{\ann}{\mathop{\rm ann}}

\providecommand{\incl}{\mathop{\rm incl}}

\newcommand{\surj}{\twoheadrightarrow}
\newcommand{\inj}{\hookrightarrow}
\newcommand{\abs}[1]{\left| #1 \right|}

\newcommand{\DS}{\displaystyle}

\begin{document}

\title[Algebraic invariants for Bestvina-Brady groups]{%
Algebraic invariants for Bestvina-Brady groups}

\author{Stefan Papadima \and Alexander I. Suciu}

\classno{%
20F36 
(primary), 
20F14,  
57M07 
(secondary)}

\extraline{The first author was partially supported by 
Program CEx 05-D11-11/2005 of the Romanian Ministry 
of Education and Research.  The second author was 
partially supported by NSF grant DMS-0311142}

\maketitle

\begin{abstract}
Bestvina-Brady groups arise as kernels of length 
homomorphisms $G_\G\to \Z$ from right-angled 
Artin groups to the integers. Under some connectivity 
assumptions on the flag complex $\Delta_\G$, we 
compute several algebraic invariants of such a group 
$N_\G$, directly from the underlying graph $\G$.  As 
an application, we give examples of finitely presented 
Bestvina-Brady groups which are not isomorphic to 
any Artin group or arrangement group.
\end{abstract}

\section{Introduction and statement of results}
\label{sect:intro}

\subsection{Bestvina-Brady groups}
\label{bbintro}
 
Given a finite simple graph $\G=(\V,\E)$, the corresponding 
right-angled Artin group $G_\G$ has presentation with 
a generator $v$ for each vertex $v\in \V$, and a 
commutator relation $vw=wv$ for each edge $\{v,w\}\in \E$.  
The Bestvina-Brady group (or, Artin kernel) associated to 
$\G$, denoted $N_{\G}$, is the kernel of the ``length" 
homomorphism to the additive group of integers, 
$\nu\colon\, G_\G \to \Z$, which sends each generator 
$v\in \V$ to $1$.  

As shown by Bestvina and Brady in their seminal 
paper \cite{BB}, the geometric and homological 
finiteness properties of the group $N_\G$ are 
intimately connected to the topology of the flag 
complex $\Delta_\G$. For example, $N_{\G}$ 
is finitely generated if and only if the graph $\G$ is 
connected; and $N_{\G}$ is finitely presented if and 
only if $\Delta_\G$ is simply-connected.  The groups 
$N_\G$ are complicated enough that a counterexample 
to either the Eilenberg-Ganea conjecture or the Whitehead 
asphericity conjecture can be constructed from them.

It is known that two right-angled Artin groups $G_{\G}$ 
and $G_{\G'}$ are isomorphic if and only if the corresponding 
graphs, $\G$ and $\G'$, are isomorphic; see \cite{KMNR}, \cite{Dr}. 
No such simple classification of the Bestvina-Brady groups 
is possible.  Indeed, if $\G$ is a tree on $n$ vertices, then 
$N_{\G}=F_{n-1}$ (the free group of rank $n-1$), 
as follows from \cite{DL}.  Thus, for any 
$n\ge 4$, there exist  graphs $\G$ and $\G'$ on $n$ 
vertices such that $\G\not\cong \G'$, yet $N_{\G}\cong N_{\G'}$.   

We study here a variety of algebraic invariants 
of a group $N_\G$ (mainly derived from the lower 
central series and the cohomology ring), showing 
how to compute these 
invariants directly from the graph $\G$, provided 
some connectivity assumptions on $\Delta_\G$ 
are satisfied. In turn, such invariants can be used to 
distinguish Bestvina-Brady groups, both among 
themselves, and from other, related classes of 
groups, such as Artin groups, or arrangement groups. 

\subsection{LCS quotients and Chen groups}
\label{lcsintro}

We start by studying  invariants derived from the 
lower central series.  For a group $G$, this series  
is defined by $\gamma_1 G=G$ and 
$\gamma_{k+1}G =(\gamma_k G,G)$, where 
$(x,y)=xyx^{-1}y^{-1}$. The direct sum of the successive 
quotients, $\gr(G)=\bigoplus\nolimits_{k\ge 1} 
\gamma_k G/ \gamma_{k+1} G$, 
is the {\em associated graded Lie algebra}\/ of $G$. 
The Lie bracket, induced from the group commutator, 
is compatible with the grading.  By construction, the 
Lie algebra $\gr(G)$ is generated by $\gr_1(G)$.
Consequently, the derived Lie subalgebra, $\gr'(G)$, 
coincides with $\bigoplus\nolimits_{k\ge 2} \gr_k(G)$. 

\begin{theorem}
\label{thm:grbbintro}
Let $\G=(\V,\E)$ be a connected graph, and let $N_{\G}$ 
be the corresponding Bestvina-Brady group.  
The associated graded Lie algebra $\gr(N_{\G})$ 
is torsion-free, with graded ranks $\phi_k=\rank \gr_k(N_{\G})$ 
given by 
\[
\prod_{k=1}^{\infty}(1-t^k)^{\phi_k}=\frac{P_{\G}(-t)}{1-t},
\]
where $P_{\G}(t)=\sum_{k\ge 0} f_k(\G) t^k$ is the clique 
polynomial of $\G$, with $f_k(\G)$ equal to the number 
of $k$-cliques of $\G$. 
Moreover, $\gr'(N_{\G})$ is isomorphic 
(as a graded Lie algebra) to the derived Lie algebra of 
$\HH_\G=\Lie( \V) / ( [v,w]=0\, \mbox{ if }\, \{v,w\}\in \E)$.
\end{theorem}

For a group $G$, let $G'=\gamma_1 G$ be the derived group, 
and $G''=(G')'$ the second derived group.  Note that 
$H_1(G)=G/G'$ is the maximal abelian quotient of $G$, 
whereas $G/G''$ is the maximal metabelian quotient.  
Define the {\em Chen Lie algebra}\/ of $G$ to be $\gr(G/G'')$.  

\begin{theorem}
\label{thm:chenbbintro}
Let $\G=(\V,\E)$ be a connected graph, and let $N_{\G}$ 
be the corresponding Bestvina-Brady group.  
The Chen Lie algebra $\gr(N_{\G}/N''_{\G})$ is 
torsion-free, with graded ranks 
$\theta_k $ given by $\theta_1=\abs{\V}-1$ and 
\[
\sum_{k=2}^{\infty} \theta_k t^{k} = 
Q_{\G} \Big(\frac{t}{1-t}\Big),
\]
where $Q_{\G}(t)=\sum_{j\ge 2} 
\left(\sum_{\W\subset \V\colon\,  \abs{\W}=j } 
\tilde{b}_0(\G_\W)\right)\! t^j$ 
is the cut polynomial of $\G$.  Moreover, 
 $\gr'(N_{\G}/N''_{\G})$ is isomorphic to the 
derived Lie algebra of  $\HH_{\G}/\HH''_{\G}$.
\end{theorem}

These two theorems rely on the analogous 
computations for right-angled Artin 
groups, done in \cite{PS-artin}. 
The proofs involve a homological analysis 
of the extension $1\to N_\G \to G_\G \to \Z \to 0$, 
based on the Salvetti complex for $G_\G$.  This 
analysis shows that $\Z$ acts trivially on $H_1(N_\G)$, 
thereby allowing us to invoke the Falk-Randell lemma \cite{FR}.

\subsection{Cohomology ring and formality}
\label{cohointro}

Next, we turn to cohomological invariants. 
If $G$ is a group, with Eilenberg-MacLane space $K(G,1)$, 
then the cohomology of $G$ with coefficients in a commutative 
ring $R$ is defined as $H^*(G,R):=H^*(K(G,1),R)$, with ring 
structure given by the cup product. The group $G$ 
is said to be $1$-formal if its Malcev Lie algebra is 
quadratically presented, cf.~\cite{Q}.  In this case, 
the rational associated graded Lie algebra 
$\gr (G)\otimes \Q$ is isomorphic 
to the rational holonomy Lie algebra, $ \HH_\Q(G)$, 
which in turn is determined by the cohomology ring 
in low degrees, $H^{\le 2}(G,\Q)$. 

For a right-angled Artin group $G_\G$, the cohomology 
ring can be identified with the exterior Stanley-Reisner 
ring of the flag complex:  $H^*(G_\G)$ is the quotient 
of the exterior algebra on generators $v^*$ in degree $1$, 
indexed by the vertices $v\in \V$, modulo the ideal generated 
by the monomials $v^* w^*$ for which $\{v,w\}$ is not an 
edge of $\G$; see \cite{KR}.  Furthermore, the group 
$G_\G$ is $1$-formal;  see \cite{KM}.

Denote by $\iota\colon\, N_\G\to G_\G$ the inclusion map 
of the kernel, and view the homomorphism 
$\nu\colon\, G_\G \to \Z\inj \Q$ as an element in 
$H^1(G_\G,\Q)$.  The next theorem determines the 
rational cohomology ring of $N_\G$, in low degrees.  
A more general result has been independently obtained 
by Leary and Saadeto\u{g}lu \cite{LeS}. 

\begin{theorem}
\label{thm:cohobbintro}
Suppose $\pi_1(\Delta_{\G})=0$. Then 
$\iota^*\colon\, H^*(G_\G,\Q)\to H^*(N_\G,\Q)$ induces 
a ring homomorphism $\iota^*\colon\, 
H^*(G_\G,\Q)/(\nu\cdot H^*(G_\G,\Q))\to H^*(N_\G,\Q)$, 
which is an isomorphism in degrees $*\le 2$.
\end{theorem}

When $\pi_1(\Delta_\G)=0$, an explicit finite presentation for 
$N_\G$ was given by Dicks and Leary \cite{DL}.  We use this 
presentation to show that the Bestvina-Brady group $N_\G$ 
is $1$-formal.

\subsection{Cohomology jumping loci}
\label{resintro}

Let $G$ be a finitely presented group, with character torus 
$\TT_G=\Hom(G,\C^*)$. Identifying the point $\rho \in \TT_G$ 
with a rank one local system ${}_{\rho} \C$ on an 
Eilenberg-MacLane space $K(G,1)$, we may define
\[ 
\label{eq=charm}
\VV(G)=\{ \rho \in \TT_G \mid  H^1(G, {} _{\rho}\C) \ne 0 \} .  
\]
The set $\VV(G)$ is an algebraic subvariety of $\TT_G$, 
called the {\em (first) characteristic variety}\/ of $G$.   
Away from the origin, this variety coincides with the 
zero set of the annihilator of the Alexander invariant, 
$B(G)\otimes \C$. 

Denote by $A=H^*(G,\C)$ the cohomology algebra of $G$. 
For each $a\in A^1$, we have $a^2=0$, and so 
right-multiplication by $a$ defines a cochain complex
$(A,a)\colon\, A^0 \stackrel{a}{\to}  A^1  \stackrel{a}{\to}   A^2$. 
Let $\RR(G)$ be the set of points $a\in A^1$ where 
this complex fails to be exact, 
\[ 
\label{eq:rv}
\RR(G)=\{a \in A^1 \mid H^1(A, a) \ne 0\}.
\]
The set $\RR(G)$ is a homogeneous algebraic 
variety in the affine space $A^1=H^1(G,\C)$, called 
the {\em (first) resonance variety}\/ of $G$.  Away from the
origin, this variety coincides with the zero set of the annihilator 
of the infinitesimal Alexander invariant, $\B(G)\otimes \C$. 

In previous work \cite{PS-artin}, \cite{DPS}, we determined the 
resonance and characteristic varieties of right-angled Artin groups. 
Here, we determine these varieties for the finitely presented 
Bestvina-Brady groups.  

Let $\TT_{\V}=(\C^*)^{\V}$ be the character torus 
of $G_\G$ (of dimension $\abs{\V}$).   For a subset 
$\W\subset \V$, let $\TT_{\W}$ be the coordinate 
subtorus supported on $\W$.  Similarly, 
let $H_{\V}=\C^{\V}$ be the Lie algebra of $\TT_\V$, 
identified with $A^1=H^1(G_\G,\C)$, and let $H_{\W}$ 
be the coordinate subspace supported on $\W$.  
The inclusion $\iota\colon\, N_{\G} \to G_{\G}$ induces 
homomorphisms 
$\iota^*\colon\, \Hom(G_{\G},\C^*) \to  \Hom(N_{\G},\C^*)$
and 
$\iota^*\colon\, H^1(G_{\G},\C) \to  H^1(N_{\G},\C)$.
Note that if $\abs{\V}=1$, then $N_{\G}=\{1\}$ and thus 
both $\VV(N_\G)$ and $\RR(N_{\G})$ are empty.  

For a graph $\G$ on vertex set $\V$, define the 
connectivity $\kappa(\G)$ to be the maximum 
integer $r$ so that, for any set of vertices $\W$ 
of size  less than $r$, the full subgraph  of $\G$ 
on vertex set $\V\setminus \W$ is connected. 

\begin{theorem}
\label{thm:resbbintro}
Let $\G$ be a graph. Suppose $\pi_1(\Delta_{\G})=0$ 
and $\abs{\V}>1$. 
\begin{enumerate}
\item
If $\kappa(\G)=1$, then $\VV(N_{\G})=\Hom(N_{\G},\C^*)$ and 
$\RR(N_{\G})=H^1(N_{\G},\C)$.  
\item
If $\kappa(\G)>1$, then the irreducible components of $\VV(N_{\G})$, 
respectively $\RR(N_{\G})$, are the subtori $\TT'_\W=\iota^*(\TT_\W)$, 
respectively the subspaces $H'_{\W}=\iota^*(H_{\W})$, 
of dimension $\abs{\W}$, one for each subset 
$\W\subset \V$, maximal among those for 
which the induced subgraph $\Gamma_{\W}$ 
is disconnected.  
\end{enumerate}
\end{theorem}

\subsection{Comparison with other classes of groups}
\label{subsec:compare}

It turns out that Bestvina-Brady groups share many common 
features with other, much-studied classes of groups:  finite-type 
Artin groups and fundamental groups of complements of 
complex hyperplane arrangements.  We catalogue here some 
of these common features, and indicate certain overlaps 
between the various classes. As a counterpoint, and as 
an application of our methods, we give examples of 
finitely presented Bestvina-Brady groups which are 
not isomorphic to any group from those two other classes.

\begin{theorem}
\label{thm:artinbbintro}
There exists an infinite family of graphs $\{\G_i\}_{i\in \N}$ 
such that the Bestvina-Brady group $N_{\G_i}$ is finitely 
presented, yet not isomorphic to either an Artin group, 
or an arrangement group. 
\end{theorem}

These graphs are obtained as the $1$-skeleta of certain 
`extra-special' triangulations of the $2$-disk. 

Using some of the machinery developed here, 
a complete classification of the Bestvina-Brady groups which 
can be realized as fundamental groups of quasi-projective 
varieties is given in \cite{DPS06}.

\subsection{Organization of the paper}
\label{org}

We start in Section \ref{sect:bb groups} with a review of 
the Bestvina-Brady groups, and a discussion of the 
Dicks-Leary presentation.  

In Section \ref{sec:salvetti} we recall the Salvetti complex 
for $G_\G$, and use it to analyze the induced homomorphism 
$\iota_*\colon\, H_1(N_\G)\to H_1(G_\G)$.

In Section \ref{sec:alex} we give presentations for the 
Alexander invariants $B(G_\G)$ and $B(N_\G)$.

In Section \ref{sect:bb lcs} we relate the graded 
Lie algebras attached to $G_\G$ and $N_\G$, 
and prove Theorems \ref{thm:grbbintro} and 
\ref{thm:chenbbintro}.

In Section \ref{sect:formal} we show that finitely presented 
Bestvina-Brady groups are $1$-formal.

In  Section \ref{sect:cohoring} we compute the homology 
groups $H_*(N_\G,\k)$, and prove Theorem \ref{thm:cohobbintro}.

In Section \ref{sect:res}  we relate the characteristic and 
resonance varieties of $G_\G$ to those of $N_\G$, and 
prove Theorem \ref{thm:resbbintro}.

Finally, in Section \ref{sect:artin-arr}, we compare finitely 
presented Bestvina-Brady groups with Artin groups and 
arrangement groups, and prove Theorem \ref{thm:artinbbintro}.

\section{Bestvina-Brady groups}
\label{sect:bb groups}

Let $\G$ be a finite graph without loops or 
multiple edges, with vertex set $\V$ and edge set 
$\E\subset {\V \choose 2}$.  The flag complex   
of $\G$, denoted $\Delta_{\G}$, is the maximal 
simplicial complex with $1$-skeleton equal to $\G$:   
the $k$-simplices of $\Delta_{\G}$ correspond 
to the $(k+1)$-cliques of $\G$. 

To the graph $\G$, there is associated a 
{\em right-angled Artin group}, $G_\G$, with a generator 
$v$ for each vertex in $\V$, and with a commutator 
relation for each edge in $\E$:  
\begin{equation}
\label{eq:artin group}
G_{\G}= \langle v \in \V \mid vw=wv 
\:  \mbox{ if }\: \{v,w\} \in \E  \rangle.
\end{equation}
For example, if $\G$ is the empty (or null) graph on $n$ vertices, 
then $G_{\G}=F_n$ (the free group of rank $n$), whereas if 
$\G$ is the complete graph $K_n$, then $G_{\G}=\Z^n$.  

\begin{definition}
\label{def:bb}
The {\em Bestvina-Brady group}\/ associated to the graph 
$\G=(\V,\E)$, denoted $N_{\G}$, is the kernel of the ``length" 
homomorphism $\nu\colon\, G_{\G}\to \Z$ which 
sends each generator $v\in \V$ to $1$.  
\end{definition}

If $\iota\colon\, N_\G \to G_\G$ denotes the inclusion map, 
we have an exact sequence
\begin{equation}
\label{eq:bb exact}
\xymatrix{1\ar[r] & N_{\G} \ar[r]^{\iota} 
& G_{\G} \ar[r]^{\nu} & \Z \ar[r] & 0}.
\end{equation}  

The group $N_{\G}$ need not be finitely generated. 
For example, if $\G$ is the empty graph on $n>1$ vertices, 
then $N_{\G}=\ker(\nu\colon\, F_n\surj \Z)$ is a free group of 
countably infinite rank.  More generally, it was shown by 
Meier--VanWyk \cite{MV} and  Bestvina--Brady \cite{BB} 
that the group $N_{\G}$ is finitely 
generated if and only if the graph $\G$ is connected. 

Even if the graph $\G$ is connected, the group $N_{\G}$ 
may not have a finite presentation. For example, if $\G$ 
is a $4$-cycle, then $G_{\G}=F_2\times F_2$; 
as noted by Stallings \cite{St}, $H_2(N_{\G})$ is not  
finitely generated, and so $N_{\G}$ is not finitely presented. 
Much more generally, Bestvina and Brady \cite{BB} 
showed that $N_{\G}$ is finitely presented if and only if 
the flag complex $\Delta_\G$ is simply-connected.  
In this case, an explicit finite presentation was given 
by Dicks and Leary \cite{DL}.   

Fix a linear order on the vertices, and orient the edges 
increasingly.  A triple of edges $(e,f,g)$ forms a directed 
triangle  if $e=\{u,v\}$, $f=\{v,w\}$, $g=\{u,w\}$, and 
$u<v<w$; see Figure \ref{fig:triangle}.

\begin{figure}
\setlength{\unitlength}{1cm}
\begin{picture}(2,1.9)(0.5,-0.4)
\xygraph{
*+{u}="u"
(
,-@{>}^{e}[ur]*+{v}="v"
,-@{>}_{g}[rr]*+{w}="w"
)
"v"-@{>}^{f}"w"
}
\end{picture}
\caption{\textsf{A directed triangle}}
\label{fig:triangle}
\end{figure}

\begin{theorem}[(Dicks--Leary \cite{DL})]
\label{thm:bb pres} 
Suppose the flag complex $\Delta_{\G}$ is simply connected. 
Then $N_\G$ has presentation 
\begin{equation}
\label{eq:bb group}
N_{\G}= \langle e \in \E \mid 
ef=fe, \: ef=g\: \mbox{ if $(e,f,g)$ is a directed triangle}\, 
\rangle.
\end{equation}
Moreover, the inclusion $\iota\colon\, N_{\G}\to G_{\G}$ 
is given by $\iota(e)=uv^{-1}$, for $e=\{u,v\}$ as above. 
\end{theorem}

The Dicks-Leary presentation is far from being 
minimal (unless $\G$ is a tree).  Indeed, there are 
$\abs{\E}$ generators in (\ref{eq:bb group}), whereas 
$H_1(N_\G)$ has rank $\abs{\V}-1$, as we shall see 
in Proposition \ref{lem=alexng}.  Nevertheless, 
(\ref{eq:bb group}) can be simplified via Tietze 
moves to a presentation where all the relations 
are commutators.  

\begin{corollary}
\label{cor:bb comm}
If $\pi_1(\Delta_\G)=0$, then $N_{\G}$ admits 
a commutator-relators presentation,  
$N_{\G}= F /R$, with $F$ the free group generated 
by the edges in a maximal tree $\mathsf{T}$, 
and $R$ a finitely generated normal subgroup of $F'$.
\end{corollary}

\begin{proof}
Fix a maximal tree $\mathsf{T}$ for $\G$.  Suppose 
$e=\{u,v\}$ is an edge not in $\mathsf{T}$.  Picking 
a path $e_1,\dots, e_r$ in $\mathsf{T}$ connecting 
$u$ to $v$, we see that 
$\iota(e)=\iota(e_1^{\epsilon_1} \cdots e_r^{\epsilon_r})$ 
in $G_\G$, for some suitable signs $\epsilon_i$. 
Thus, $e=e_1^{\epsilon_1} \cdots e_r^{\epsilon_r}$ in $N_\G$.  
This shows that $N_\G$ is generated by the edges 
of $\mathsf{T}$.  
Now note that $N_\G/N'_\G$ is free abelian, of rank 
equal to the number of edges in 
$\mathsf{T}$.  Eliminating the redundant generators 
from (\ref{eq:bb group}), we arrive at the desired 
presentation.   
\end{proof}

In certain situations, the Dicks-Leary presentation 
permits us to identify the group $N_\G$ in terms 
of better known groups. 

\begin{example}
\label{ex:bb tree}
Suppose $\G$ is a tree on $n$ vertices.  Then $\G$ has 
no triangles, and $\pi_1(\Delta_\G)=0$.   Since 
$\Gamma$ has $n-1$ edges, we see that $N_{\G}=F_{n-1}$.  
\end{example}

\begin{example}
\label{ex:bb cone}
Suppose $\G$ is the cone on $\G'$.  Then 
$G_\G=G_{\G'}\times \Z$, and so $N_\G=G_{\G'}$. 
In particular, if $\G=K_n$, then $N_\G=\Z^{n-1}$.
\end{example}

In general, though, $N_{\G}$ is not isomorphic to 
any right-angled Artin group, as we shall show later. 

Noteworthy is the situation when $\Delta_\G$ is 
a triangulation of the $2$-disk.  In this case, $N_\G$ 
admits a $2$-dimensional $K(N_\G,1)$, see \cite[Corollary 2.3]{Cr}.   

\begin{definition}
\label{def:special triang}
A triangulation of the disk is said to be {\em special}\/
if it is obtained from a triangle by adding one triangle 
at a time, along a unique boundary edge. 
\end{definition}

\begin{lemma}
 \label{lem:special triang}
 Let $\Delta$ be a special triangulation of $D^2$, 
with $1$-skeleton $\G=(\V,\E)$. Then:
 \begin{enumerate}
 \item \label{sp1}
 $2\abs{\V}-\abs{\E}=3$.
\item \label{sp2} 
$\Delta_\G=\Delta$.
 \item  \label{sp3}
 $N_\G$ admits a presentation with 
 $\abs{\V}-1$ generators and $\abs{\V}-2$ 
 commutator relators.
  \end{enumerate}
 \end{lemma}
 
 \begin{proof}
By induction on the number $t$ of triangles.  
Evidently, all statements hold for a single triangle.  
Now suppose $\Delta$ is a special triangulation 
with $t$ triangles, and a directed triangle $(e,f,g)$ 
is added along edge $e$, to form $\Delta'$. 
In the process, one vertex and two edges are added, 
and so  the quantity $2\abs{\V}-\abs{\E}$ does not 
change.  The only new $3$-cycle in the graph $\Gamma'$ 
is the boundary of  $(e,f,g)$; thus, $\Delta'$ is a flag complex. 
Furthermore, if $\T$ is a maximal tree for $\Gamma$, we can 
build a maximal tree $\T'$ for $\Gamma'$ by adding a  
new edge, say, $f$. The Dicks-Leary relation $ef=g$ 
(with $e$ expressed as a word in the edges of $\T$) 
may be used to eliminate the generator $g$.  Thus, 
 $N_{\G'}$ has only one new generator, $f$, and only 
 one new relation, $ef=fe$.  
 \end{proof}

\begin{figure}
\setlength{\unitlength}{1cm}
\begin{picture}(4,3.1)(0.8,-0.4)
\xygraph{
*+{v_1}="v1"
(
,-@{>}[ur]*+{v_4}="v4"
,-@{>}_{e_1}[rr]*+{v_2}="v2"
(
,-@{>}_{e_4}[rr]*+{v_3}="v3"
,-@{>}_{e_3}[ur]*+{v_5}="v5"-@{>}_{e_5}[ul]*+{v_6}="v6"
)
)
"v2"-@{>}^{e_2}"v4"
"v3"-@{>}"v5"
"v4"-@{>}"v5"
"v4"-@{>}"v6"
}
\end{picture}
\caption{\textsf{A special triangulation of the $2$-disk}}
\label{fig:big triangle}
\end{figure}

\begin{example}
\label{ex:double triangle}
Let $\G$ be the graph in Figure \ref{fig:big triangle}. 
Choosing a maximal tree $\T=\{e_1,\dots ,e_5\}$ as 
indicated, the presentation from 
Lemma \ref{lem:special triang}\ref{sp3} 
reads as follows:
\[
N_{\G}=\langle e_1,\dots ,e_ 5 \mid 
(e_1,e_2), \ (e_2,e_3), \ (e_3,e_4),\ 
e_5e_2^{-1}e_3=e_2^{-1}e_3e_5\rangle.
\]
We shall see in Proposition \ref{prop:special bb artin}  
that $N_{\G} \not\cong G_{\G'}$, no matter what the graph $\G'$ is. 
\end{example}

\section{The Salvetti complex}
\label{sec:salvetti}

For a simple graph $\G=(\V, \E)$, let $K_{\G}$ be the 
CW-complex obtained by joining tori in the manner 
prescribed by the flag complex $\Delta_{\G}$.  
More precisely, if $T^n=(S^1)^{\times n}$ is the torus of 
dimension $n=\abs{\V}$, with the usual CW-decompo\-sition, 
then $K_{\G}$ is the subcomplex obtained by deleting 
the cells corresponding to the non-faces of $\Delta_{\G}$. 

Clearly, the fundamental group of $K_{\G}$ 
is the right-angled Artin group $G_\G$. 
In fact, $K_{\G}$ is an Eilenberg-MacLane space of type 
$K(G_{\G},1)$; see \cite{MV}.  Note that $H_k(K_{\G})$ 
is free abelian, of rank equal to the number of $k$-cliques 
in $\G$; thus, the Poincar\'e polynomial of $K_{\G}$ 
equals the clique polynomial of the graph, $P_{\G}(t)$. 

We use the cell structure of the space $K_\G$ to describe 
a finite, free resolution of $\Z$, viewed as a 
trivial module over the group ring $\Z G_{\G}$. This resolution 
was first determined by Salvetti \cite{Sa}, in the more general 
context of Artin groups, and further extended by Charney 
and Davis \cite{CD}.  For the benefit of the reader, we 
include a self-contained, direct computation of the 
Salvetti complex in our particular case.

\subsection{A free $\Z G_{\G}$-resolution of $\Z$}
\label{subsec:free res} 
For each subset $\W\subset \V$, let 
$\G_{\W}$ be the induced subgraph of $\G$ on vertex set $\W$.   
Let $G_{\W}=G_{\G_{\W}}$ be the corresponding 
right-angled Artin group, and let $K_{\W}=K_{\G_{\W}}$ 
be the corresponding CW-complex.  
The inclusion $\W\subset \V$ gives rise to a cellular 
inclusion map $j_{\W}\colon\, K_{\W}\to K_{\G}$.  
The induced homomorphism, 
$(j_{\W})_{\#} \colon\, G_{\W}\to G_{\G}$,   
is a split injection, with retract 
 $G_{\G}\to G_{\W}$ given on generators 
 by $v\mapsto v$ if $v\in \W$, and $v\mapsto 1$ otherwise.

Denote by $C_{\bullet}=C_{\bullet}(K_{\G})$ the 
cellular chain complex of  $K_{\G}$.  A basis for 
$C_k$ is given by the complete $k$-subgraphs of $\G$:  
to each complete subgraph on vertex set $\W\subset \V$, 
there corresponds a cell $c_{\W}$. Since $C_{\bullet}$ 
is a sub-complex of the cellular chain complex of $T^n$, 
the boundary maps $C_k\to C_{k-1}$ are the zero maps.  

Now let $\widetilde{C}_{\bullet}=
(C_{\bullet}(\widetilde{K}_{\G}),\partial_{\bullet})$ be the 
equivariant chain complex of the universal cover of $K_{\G}$.  
The augmentation map, $\epsilon\colon\, \widetilde{C}_0 \to \Z$, 
extends to a finite, free resolution $\widetilde{C}_{\bullet} \to \Z$ 
over the group ring $\Z G_{\G}$.  

\begin{prop}
\label{prop:chain complex}
Under the identification $\widetilde{C}_k = \Z G_{\G} \otimes C_k$, 
the boundary map $\partial_k\colon\, 
\widetilde{C}_k\to \widetilde{C}_{k-1}$ is given by:
\begin{equation}
\label{eq:bdry map}
\partial_k (c_{\W}) = \sum_{r=1}^{k}  
(-1)^{r-1} (v_{i_r} -1) c_{\W\setminus \{i_r\} }
\end{equation}
where $\W=\{v_{i_1},\dots , v_{i_k}\}$ is a $k$-clique in $\G$. 
\end{prop}

\begin{proof}
Let $K_{\W}=T^k$ be the corresponding CW-subcomplex 
of $K_{\G}$. The equivariant chain complex 
$\widetilde{C}^{\W}_{\bullet}=
(C_{\bullet}(\widetilde{K}_{\W}),\partial^{\W}_{\bullet})$ 
is simply the Koszul complex on the variables from $\W$.  In particular, 
 $\partial^{\W}_k (c_{\W})$ is given by the right-hand side of 
(\ref{eq:bdry map}). The commutativity of the diagram
\[
\xymatrix{
\Z G_{\W}  \otimes C_k(K_{\W}) \ar[r]^{\partial^{\W}_k} 
\ar[d]^{j_{\#} \otimes j_*} & 
\Z G_{\W}  \otimes C_{k-1}(K_{\W}) 
\ar[d]^{j_{\#} \otimes j_*} \\
\Z G_{\G}  \otimes C_k(K_{\G}) \ar[r]^{\partial_k} & 
\Z G_{\G}  \otimes C_{k-1}(K_{\G}) 
}
\]
completes the proof.
\end{proof}

\subsection{An injectivity lemma}
\label{subsec:action} 
Recall that the Bestvina-Brady group associated to a graph 
$\G=(\V,\E)$ is the kernel of the homomorphism 
$\nu\colon\, G_{\G}\to \Z$ that sends each generator $v\in \V$ 
of $G_{\G}$ to $1\in \Z$.  Let $\iota\colon\, N_\G \to G_\G$ 
be the inclusion map.   

\begin{lemma}
\label{lem:inj map}
If the graph $\G$ is connected, then the induced 
homomorphism $\iota_*\colon\, H_1(N_{\G}) \to H_1(G_{\G})$ 
is injective.
\end{lemma}

\begin{proof}
Let $\widetilde{C}_{\bullet}=
(C_{\bullet}(\widetilde{K}_{\G}),\partial_{\bullet})$ 
be the Salvetti complex, with boundary maps given 
by (\ref{eq:bdry map}).  Identify the group ring 
$\Z\Z$ with the ring of Laurent polynomials 
$\Z[\tau,\tau^{-1}]$.  By Shapiro's Lemma, 
$H_*(N_{\G})=H_*(\Z\Z\otimes _{\Z G_{\G}} 
\widetilde{C}_{\bullet} )$, where 
$\Z\Z=\Z[\tau,\tau^{-1}]$ is viewed as a right 
$G_{\G}$-module via the map $v\mapsto \tau$. 

After identifying $\Z\Z\otimes _{\Z G_{\G}} \widetilde{C}_k$ 
with $\Z\Z\otimes C_k$, the boundary map $\partial_k$ 
takes the form $(\tau-1)\otimes d_k\colon\, 
\Z\Z\otimes C_k \to \Z\Z\otimes C_{k-1}$, 
where $d_k$ is the simplicial boundary 
map on the $(k-1)$-simplices of $\Delta_{\G}$. 
With these identifications,  the chain map 
$\iota_{\#}\colon\, \Z\Z\otimes_{\Z G_{\G}} 
C_{\bullet}(\widetilde{K}_{\G})\to  C_{\bullet}(K_{\G})$ 
takes the form:
\[
\xymatrixcolsep{18pt}
\xymatrix{ 
\dots \ar[r]&\Z \Z \otimes \Z^{\E} 
\ar[rrr]^{\partial_2=(\tau-1) \otimes d_2} 
\ar[d]^{\iota_2=\epsilon \otimes \id} &&&
\Z \Z  \otimes \Z^{\V}
\ar[rrr]^{\partial_1=(\tau-1) \otimes d_1} 
\ar[d]^{\iota_1=\epsilon \otimes \id} &&&
\Z \Z  \otimes \Z \ar[d]^{\iota_0=\epsilon \otimes \id}
\\
\cdots \ar[r]&\Z   \otimes \Z^{\E} \ar[rrr]^{0} &&& 
\Z  \otimes \Z^{\V} \ar[rrr]^{0}  &&& 
\Z  \otimes \Z
}
\]

The map $\iota_*\colon\, H_1(N_{\G}) \to H_1(G_{\G})$ is the 
homomorphism induced on $H_1$ by the middle down  
arrow. To show $\iota_*$ is injective, we need to prove: 
$\ker \partial_1\cap \ker \iota_1 \subset \im \partial_2$. 

Let $z=\sum_{v\in \V} p_v \otimes v \in \Z\Z \otimes \Z^{\V}$.   
Suppose $z$ belongs to $\ker \partial_1$.  
Since $\partial_1(z)=(\tau-1) (\sum_{v\in \V} p_v)$, and since 
the ring $\Z\Z$ has no zero-divisors, we must have 
$\sum_{v\in \V} p_v=0$. 
Now suppose $z$ belongs to $\ker \iota_1$.  
Then $\sum_{v\in \V} \epsilon(p_v)  \otimes v=0$, which can only 
happen if $\epsilon(p_v)=0$, for all $v\in \V$.  Thus, for each 
$v\in \V $, there is $q_v\in \Z\Z$ such that $p_v=(\tau-1) q_v$. 
We conclude that $\ker \partial_1\cap \ker \iota_1$ 
is generated by elements of the form 
$(\tau-1)q\otimes (v- u)$.

Let $e_1,\dots ,e_s$ be a path in $\G$ 
with $e_i=(u_{i-1},u_i)$, joining $u_0=u$ to $u_s=v$.  
Then:
\[
\partial_2 \Big(\sum_{i=1}^{s} q\otimes e_i \Big) = 
(\tau-1) \sum_{i=1}^{s} q\otimes (u_i-u_{i-1}) = 
(\tau-1) q \otimes (v-u).
\]
This finishes the proof.
\end{proof}

\section{Alexander invariants}
\label{sec:alex}

Let $G$ be a group, with abelianization $H_1(G)=G/G'$.  
The {\em Alexander invariant}\/ of $G$ is the quotient 
group $B(G)= G'/G''$, endowed with the $\Z H_1(G)$-module 
structure induced by conjugation in $G/G''$, via the exact 
sequence $0\to G'/G''\to G/G''\to G/G'\to 0$.  Alternatively, 
if $K$ is a connected CW-complex with $\pi_1(K)=G$, 
and $K^{\ab}$ is the universal abelian cover of $K$, 
then $B(G)=H_1(K^{\ab})$, with module structure 
coming from the action of $H_1(G)$ by deck 
transformations. 

In this Section, we determine the Alexander 
invariants of right-angled Artin groups and of finitely-generated 
Bestvina-Brady groups. We start by giving a finite presentation 
for $B({G_{\G}})$, viewed as a module over $\Z H_1(G_{\G})$, 
for an arbitrary finite graph $\G=(\V, \E)$. 

After fixing a total ordering on $\V$, we may 
identify $\Z H_1(G_{\G})$ with the ring of Laurent 
polynomials in variables labeled by the vertices,  
$\Lambda=\Z[\V^{\pm 1}]$. Let $t$ be a triple of vertices, 
and let $e$ be a $2$-element subset of $t$. 
We denote by $v_e$ the third vertex of $t$, 
and by $\epsilon_e$ the sign of the permutation $(v_e, u, w)$ 
where $e=\{u,w\}$, with $u<w$. 

\begin{theorem}
\label{thm=alexpres}
The Alexander invariant $B(G_{\G})$ is the 
$\Z[\V^{\pm 1}]$-module generated by the 
non-edges $e\in \E_{\overline{\G}}$, and with relators 
\[
\sum_{e\subset t, \, e\in \E_{\overline{\G}}}
\epsilon_e (v_e-1) \otimes e, 
\]
indexed by the triples of vertices $t\in {\V \choose 3}$ which are not 
$3$-cliques of $\G$.
\end{theorem}

\begin{proof}
As before, let $\widetilde{C}_{\bullet}=
(C_{\bullet}(\widetilde{K}_{\G}),\partial_{\bullet})$
be the equivariant chain complex of $K_{\G}$.
By Shapiro's Lemma, $B({G_\G})=
H_1(\Z H_1(G_\G) \otimes_{\Z G_\G} \widetilde{C}_{\bullet})$, 
as $\Lambda$-modules. 

Recall $K_{\G}$ is a CW-subcomplex of the torus 
$T^n$, where $n=\abs{\V}$.  Since $T^n=K(\Z^n,1)$, the 
equivariant chain complex 
$(C_{\bullet}(\widetilde{T}^n),\delta_{\bullet})$ gives   
a free $\Lambda$-resolution of $\Z$.  Using 
(\ref{eq:bdry map}), it is readily seen that  $\Lambda
\otimes_{\Z G_\G} \widetilde{C}_{\bullet} =
\Lambda \otimes C_{\bullet}$ is a 
$\Lambda$-subcomplex of 
$C_{\bullet}(\widetilde{T}^n) = 
\Lambda^{\V \choose \bullet}$, and that these two 
complexes coincide up to degree $\bullet=1$.  
A diagram chase yields 
\begin{equation}
\label{eq=unred}
B({G_\G})= \coker \left(\delta_3 + \incl \colon\,  
\Lambda^{\V \choose 3} \oplus  \Lambda^{\E} 
 \to\Lambda^{\V \choose 2} \right).
\end{equation}
Row-reducing the above presentation matrix finishes the proof. 
\end{proof}

A homomorphism  $\iota\colon\, N \to G$ induces in a natural 
way a change of rings map $\iota_*\colon\, \Z H_1(N)\to 
\Z H_1(G)$ and a $\Z H_1(N)$-linear map $B(\iota) \colon\, 
B(N) \to B(G)$.

\begin{prop}
\label{lem=alexng}
Let $\G$ be a connected graph, and let 
$\iota\colon\, N_\G\to G_\G$ be the inclusion map.
Then,  the following hold.
\begin{enumerate}
\item \label{ng1}
There is a split exact sequence
$\xymatrixcolsep{16pt}
\xymatrix{0\ar[r] & H_1(N_{\G}) \ar[r]^{\iota_*} 
& H_1(G_{\G}) \ar[r]^(.6){\nu_*} & \Z \ar[r] & 0 }$. 

\item \label{ng0}
The inclusion $\iota$ restricts to an equality 
$N_{\G}'=G'_\G$.

\item \label{ng2}
The induced map $B(\iota) \colon\, 
B(N_{\G}) \to B(G_{\G})$ is a $\Z H_1(N_\G)$-linear 
isomorphism. 
\end{enumerate}
\end{prop}

\begin{proof}
(\ref{ng1})   This a direct consequence of 
Lemma \ref{lem:inj map}. 

(\ref{ng0})  Clearly, $N_{\G}' \subset G'_\G$.   
Now let $g\in G'_\G$; then $\nu(g)=0$, and so $g=\iota(n)$, for 
some $n\in N_{\G}$. Since $\iota_*\colon\, H_1(N_{\G})\to H_1(G_\G)$ 
is injective, we must have $n\in N_{\G}'$. 

(\ref{ng2}) 
By the above, $N_{\G}''=G_{\G}''$, and so 
 $N_{\G}'/N_{\G}''=G_{\G}'/G_{\G}''$.
\end{proof}

\begin{corollary}
\label{cor=nalex}
Let $\G$ be a connected graph, with vertex set $\V=\{1,\dots,n\}$. 
The Alexander invariant $B(N_{\G})$ is isomorphic to the 
restriction of the $\Z\Z^n$-module $B(G_{\G})$, with 
presentation given in Theorem \ref{thm=alexpres}, 
via the change of rings $\iota_*\colon\, \Z\Z^{n-1}\to \Z\Z^n$. 
\end{corollary}

Let $G$ be a finitely presented group, with torsion-free 
abelianization.  
The  {\em holonomy Lie algebra}\/ of $G$, 
denoted $\HH(G)$, is the quotient of the free Lie algebra 
on $H_1(G)$ by the ideal generated by the image of the 
comultiplication map 
$\nabla_G\colon\, H_2(G) \to H_1(G)\wedge H_1(G)$. 

The {\em infinitesimal Alexander 
invariant}\/ of $G$ is $\B(G)=\HH(G)'/\HH(G)''$, 
with module structure over the symmetric 
algebra $\Sym(H_1(G))$ coming from 
the exact sequence 
\[
0\to \HH(G)'/\HH(G)''\to \HH(G)/\HH(G)''\to 
\HH(G)/\HH(G)'\to 0.
\]  
This module is isomorphic to the ``linearization" of the 
classical Alexander invariant of the group $G$, see \cite{PS-chen}. 

A homomorphism $\iota\colon\, N\to G$ induces 
 a change of rings map $\iota_*\colon\, \Sym(H_1(N))\to 
 \Sym(H_1(G))$ and a $ \Sym(H_1(N))$-linear map 
 $\B(\iota)\colon\, \B(N)\to \B(G)$.

Now let $\G=(\V,\E)$ be a finite graph.  After fixing a total 
ordering on $\V$, we may identify $\Sym(H_1(G_\G))$ with 
the polynomial ring $S=\Z[\V]$.   
Using \cite[Theorem 6.2]{PS-chen}, we obtain the following 
infinitesimal analogue of Theorem \ref{thm=alexpres}.

\begin{prop}
\label{prop:clex inv artin}
The infinitesimal  Alexander invariant of a right-angled 
Artin group, $\B(G_{\G})=\HH(G_{\G})'/\HH(G_{\G})''$, 
is the $S$-module generated by the  non-edges 
$e\in \E_{\overline{\G}}$, and with relators 
\[
\sum_{e\subset t, \, e\in \E_{\overline{\G}}} 
\epsilon_e v_e \otimes e, 
\]
indexed by the triples $t\in {\V\choose 3}$ which are not 
$3$-cliques of $\G$. 
\end{prop}

\section{Lower central series}
\label{sect:bb lcs}

In this section, we determine the associated graded Lie algebra 
and the Chen Lie algebra of the Bestvina-Brady group 
corresponding to a finite, connected graph, thus proving 
Theorems \ref{thm:grbbintro} and \ref{thm:chenbbintro} 
from the Introduction. 

\subsection{Lie algebras associated to right-angled Artin groups}
\label{subsec:artin lcs}

A graph $\G=(\V,\E)$ determines in a natural way 
a graded, finitely-presented Lie algebra $\HH_\G$, 
as follows: 
\begin{equation}
\label{eq:holo artin}
 \HH_\G=\Lie( \V) / ( [v,w]=0\ \mbox{if $\{v,w\}\in \E$}), 
\end{equation}
where $\Lie( \V)$ is the free Lie algebra on the vertex set  $\V$. 

For each $k\ge 1$, let $f_k(\G)$ be the number of complete 
$k$-subgraphs of $\G$, and set $f_0(\G)=1$. 

\begin{theorem}[(\cite{DK1}, \cite{DK2}, \cite{PS-artin})]
\label{thm:lcs artin}
Let $\G$ be a finite graph, and let 
$G_\G$ be the corresponding right-angled 
Artin group. Then  $\gr(G_{\G})\cong \HH_\G$, 
as graded Lie algebras. 
Moreover,  the graded pieces of $\gr(G_{\G})$ 
are torsion-free, with ranks $\phi_k=\rank (\gr_k(G_{\G}))$  
given by:
\[
\label{eq:lcs artin}
\prod_{k=1}^{\infty}(1-t^k)^{\phi_k}=P_{\G}(-t),
\]
where $P_{\G}(t)=\sum_{k\ge 0} f_k(\G) t^k$ is the clique 
polynomial of  $\G$. 
\end{theorem}

For each $j\ge 1$, let 
$c_j(\G)=\sum_{\W\subset \V\colon\,  \abs{\W}=j } \tilde{b}_0(\G_\W)$, 
where $\tilde{b}_0(\G)=\rank \widetilde{H}_0(\G)$ is 
the number of components of $\G$ minus $1$.  
Note that $c_1(\G)=0$, and also $c_j(\G)=0$, if 
$j>  \abs{\V}-\kappa(\G)$, where $\kappa(\G)$ is the 
connectivity of $\G$. 

\begin{theorem}[\cite{PS-artin}] 
\label{thm:chen artin}
Let $\G$ be a finite graph. 
Then $\gr(G_{\G}/G_{\G}'')\cong \HH_\G/\HH_\G''$, 
as graded Lie algebras.  Moreover,  the graded pieces of 
$\gr(G_{\G}/G_{\G}'')$ are torsion-free, with ranks 
$\theta_k=\rank (\gr_k(G_{\G}/G_{\G}''))$  given by:
\[
\label{eq:chen ranks artin}
\sum_{k=2}^{\infty} \theta_k t^{k} = 
Q_{\G} \Big(\frac{t}{1-t}\Big),
\]
where $Q_{\G}(t)=\sum_{j= 2}^{\abs{\V}-\kappa(\G)} c_j(\G) t^j$ 
is the cut polynomial of $\G$.
\end{theorem}

\subsection{Monodromy action}
\label{subsec:mono action}
Let $N_\G$ be the Bestvina-Brady group associated 
to the graph $\G$.  Recall we have an exact sequence 
\begin{equation}
\label{eq:bb}
\xymatrix{1\ar[r] & N_{\G} \ar[r]^{\iota} 
& G_{\G} \ar[r]^{\nu} & \Z \ar[r] & 0}.
\end{equation}  
This sequence admits a splitting $s \colon\, \Z\to G_{\G}$, 
given by $s(1)=v$, for some fixed generator $v\in \V$.  

\begin{prop}
\label{prop:trivial action}
Let $\G$ be a connected graph.  Then, in the split extension 
(\ref{eq:bb}), the group $\Z$ acts trivially on the abelianization 
$H_1(N_{\G})$.  
\end{prop}

\begin{proof}
The monodromy of the semidirect product 
$G_{\G}=N_{\G}\rtimes_{\sigma} \Z$ is given by 
$\sigma(1)(x)=vxv^{-1}$, for some fixed generator 
$v\in \V$. Conjugation by any element of 
$G_{\G}$ acts trivially on $H_1(G_{\G})$.  On the 
other hand, we know from Lemma \ref{lem:inj map}
that $H_1(N_{\G})$ injects into $H_1(G_{\G})$.  Hence,  
conjugation by an element of $G_\G$ also acts 
trivially on $H_1(N_{\G})$. 
\end{proof}

For a homomorphism $\alpha\colon\, G\to H$, let 
$\bar\alpha\colon\, G/G'' \to H/H''$ be the induced 
homomorphism on maximal metabelian quotients.

\begin{prop}
\label{prop:gr chen bb}
Let $\G$ be a connected graph.  Then,  the sequence
  \begin{equation}
\label{eq:bb chen}
\xymatrix{1\ar[r] & N_{\G}/N_{\G}'' \ar[r]^{\bar\iota} 
& G_{\G}/G_{\G}'' \ar[r]^(.6){\bar\nu} & \Z \ar[r] & 0 },
\end{equation}
is split exact, with trivial monodromy action 
on $H_1(N_{\G}/N_{\G}'')=H_1(N_{\G})$.
\end{prop}

\begin{proof}
From the proof of Proposition \ref{lem=alexng}, 
we know that $N_{\G}''=G_{\G}''$. 
All claimed properties of sequence (\ref{eq:bb chen}) 
follow from the corresponding properties of (\ref{eq:bb}). 
\end{proof}

\subsection{Lie algebras associated to Bestvina-Brady groups}
\label{subsec:bb lcs}

Recall now the following well-known result of 
Falk and Randell. 

\begin{theorem}[(\cite{FR})]
\label{fr}
Let $\xymatrixcolsep{15pt}
\xymatrix{1\ar[r] &  A \ar^{\alpha}[r] 
& B \ar^{\beta}[r] & C \ar[r] & 1}$ 
be a split exact sequence of groups.  
Suppose $C$ acts trivially on $H_1(A)$. Then 
\[
\xymatrix{0\ar[r] &  \gr(A) \ar^{\gr(\alpha)}[r] &
\gr(B) \ar^{\gr(\beta)}[r] & \gr(C) \ar[r] & 0}
\]
is a split exact sequence of graded Lie algebras. 
\end{theorem}

This result permits us to reduce the computation of the 
LCS and Chen quotients of a Bestvina-Brady group 
to the computation of the LCS and Chen quotients 
of the corresponding right-angled Artin group. 

\begin{theorem}
\label{thm:gr chen bb}
Let $\G$ be a finite, connected graph.  Then, 
the inclusion map $\iota\colon\, N_{\G}\to G_{\G}$ induces 
isomorphisms of graded Lie algebras
\begin{enumerate}
\item \label{bb1}
$\gr'(\iota) \colon\, \gr'(N_{\G}) \stackrel{\,\cong\,}{\longrightarrow}
\gr'(G_{\G})$. 
\item \label{bb2}
$\gr'(\bar\iota) \colon\, \gr'(N_{\G}/N''_{\G}) 
\stackrel{\,\cong\,}{\longrightarrow} \gr'(G_{\G}/G''_{\G})$. 
\end{enumerate}
\end{theorem}

\begin{proof}
Using Propositions \ref{prop:trivial action} and 
\ref{prop:gr chen bb}, we may apply Theorem \ref{fr} to 
the exact sequences (\ref{eq:bb}) and (\ref{eq:bb chen}).  
Noting that $\gr(\Z)=\Z$ (concentrated in degree $1$), 
yields isomorphisms \ref{bb1} and \ref{bb2}, 
respectively. 
\end{proof}

Combining Theorem \ref{thm:gr chen bb} with 
Theorems \ref{thm:lcs artin} and \ref{thm:chen artin} 
finishes the proof of Theorems \ref{thm:grbbintro} 
and \ref{thm:chenbbintro} from the Introduction.

\begin{example}
\label{ex:bb tree 2}
Suppose $\G$ is a tree on $n$ vertices.  Then recall 
$N_{\G}$ is a free group of rank $n-1$.  
By Theorem \ref{thm:gr chen bb}, 
$\phi_k(G_{\G})=\phi_k( F_{n-1})$ and 
$\theta_k(G_{\G})=\theta_k( F_{n-1})$, for 
all $k\ge 2$, which recovers the computations from 
\cite[\S 6.2]{PS-artin}. 
\end{example}

\begin{remark}
\label{rem:bb}
Suppose  $\G=K_{n_1,\dots ,n_r}$ is a complete multi-partite 
graph. Then $G_{\G}=F_{n_1}\times \cdots \times F_{n_r}$, 
and so, by Theorem  \ref{thm:gr chen bb}, $\gr' (N_{\G})
= \gr'(F_{n_1}\times \cdots \times F_{n_r})$.  
Even though, from the point of view 
of the lower central series quotients, $N_{\G}$ looks 
like a product of free groups, it is not of this type, 
except when some $n_i=1$.  Indeed, if all $n_i>1$, 
then $H_{r-1}(\Delta_{\G})$ is a free abelian group of rank 
$\prod_{i=1}^{r} (n_i-1)>0 $, and so it follows from \cite{BB} 
that $N_{\G}$ does not have a finite $K(N_{\G},1)$.  
\end{remark}

\section{Holonomy Lie algebra and $1$-formality}
\label{sect:formal}

Let $G$ be a finitely presented group, with torsion-free 
abelianization.    Recall that the holonomy Lie algebra 
$\HH(G)$ is the quotient of the free Lie algebra 
$\Lie(H_1(G))$ by the ideal generated by the image of 
the comultiplication map 
$\nabla_G\colon\, H_2(G) \to H_1(G)\wedge H_1(G)$. 
Note that $\HH(G)$ inherits a natural grading from the 
free Lie algebra, compatible with the Lie bracket.  
By construction, $\HH(G)$ is generated by $\HH_1(G)$.
Consequently, the derived Lie subalgebra, $\HH'(G)$, 
coincides with $\HH_{\ge 2}(G)$. 
If we drop the torsion-freeness assumption on $H_1(G)$, 
we may still define the rational holonomy Lie algebra, 
$\HH_\Q(G)$, using the rational homology groups of $G$. 

To a group $G$, Quillen associates in a functorial way 
a Malcev filtered Lie algebra, $M_G$; 
see \cite[Appendix A]{Q}, 
and also \cite{PS-chen} for further details.  
A finitely presented group $G$ is said to be 
{\em $1$-formal}\/ if $M_G$ is isomorphic to the 
rational holonomy Lie algebra, $\HH_\Q(G)$, 
completed with respect to the bracket length filtration. 
Equivalently, $M_G$ is a quadratic Malcev Lie algebra.
If the group $G$ is $1$-formal, then
$\gr (G)\otimes \Q\cong \HH_\Q(G)$, 
as follows from \cite{Q}.  Moreover, 
$\gr (G/G'')\otimes \Q\cong \HH_\Q(G)/\HH_\Q(G)''$, 
as shown in \cite{PS-chen}.  

Assume now $G$ is finitely presented, and $H_1(G)$ is torsion-free.  
Then, the canonical projection $\Lie(H_1(G)) \surj \gr(G)$ 
factors through an epimorphism of graded Lie algebras, 
$\Psi_G \colon\,  \HH(G) \surj \gr (G)$, which  in turn  
descends to an epimorphism 
\[
\Psi^{(2)}_G \colon\,  \HH(G)/\HH(G)''  \surj \gr (G/G'').
\]  
If the group $G$ is $1$-formal, then the maps 
$\Psi_G\otimes \Q$ and $\Psi^{(2)}_G\otimes \Q$ 
are isomorphisms, see \cite{PS-chen}. 

By construction, the resonance variety $\RR(G)$ depends 
only on the holonomy Lie algebra $\HH_\Q(G)$. More precisely, 
if  $\HH_\Q(G_1)\cong \HH_\Q(G_2)$, as graded Lie algebras, 
then there is a linear isomorphism $H^1(G_1,\C) \cong H^1(G_2,\C)$, 
restricting to an isomorphism $\RR(G_1)\cong \RR(G_2)$.  

For a right-angled Artin group $G_\G$, it is easily seen 
that $\HH(G_\G)=\HH_\G$, cf.~\cite{PS-artin}. Moreover, 
as shown by Kapovich and Millson \cite{KM}, the group 
$G_{\G}$ is $1$-formal. We now prove an analogous 
result for the Bestvina-Brady groups. 

\begin{prop}
\label{prop:bb formal}
If the flag complex $\Delta_{\G}$ is simply-connected,  
then $N_{\G}$ is $1$-formal.
\end{prop}

\begin{proof}
Consider the Dicks-Leary  presentation (\ref{eq:bb group}) 
for $N_\G$. It follows from \cite{P1} that the Malcev Lie 
algebra of $N_{\G}$ is the quotient of  $\widehat{\Lie(\E)}$, 
the free Malcev Lie algebra on $\E$,
 by the closed Lie ideal generated by 
the elements $(e,f), (f,g), (e, g), efg^{-1}$, for all 
directed triangles $(e,f,g)$ as in Figure \ref{fig:triangle}. 
Here multiplication denotes 
the Campbell-Hausdorff product in the underlying 
Malcev group.  

Now use  \cite[Lemma 2.5]{P2} to replace the CH 
commutators $(e,f), (f,g), (e, g)$ by the corresponding 
Lie commutators, $[e,f], [f,g], [e, g]$.  It follows from 
the definition of the CH product that we may also 
replace $efg^{-1}$ by $e+f-g$. 
This shows that $M_{N_{\G}}$ is a quadratic 
Malcev Lie algebra, and so, $N_{\G}$ is $1$-formal.
\end{proof}

Let  $\iota\colon\, N \to G$ be a homomorphism between 
finitely presented groups with torsion-free abelianizations.  
Denote by $\iota_*\colon\, H_*(N) \to H_*(G)$ the induced 
homomorphism in homology.  We then  have a commuting 
diagram, 
 \begin{equation}
 \label{eq:nabla}
 \xymatrix{
 H_2(N) \ar[r]^(.35){\nabla_N} \ar[d]^{\iota_*} & 
  H_1(N)\wedge H_1(N) \ar[d]^{\iota_*\wedge \iota_*}\\
   H_2(G) \ar[r]^(.35){\nabla_G} &
  H_1(G)\wedge H_1(G)
 }
\end{equation}
Consequently, there is an induced morphism of 
graded Lie algebras, 
$\HH(\iota) \colon\,  \HH(N) \to \HH(G)$.

\begin{lemma}
\label{lem:hhprime}
If $\pi_1(\Delta_{\G})=0$, then the map 
$\HH'_\Q(\iota)\colon\, \HH'_\Q(N_\G) \to \HH'_\Q(G_\G)$ 
is an isomorphism of graded Lie algebras. 
\end{lemma}

\begin{proof}
The inclusion $\iota\colon\, N_\G\to G_\G$ induces Lie algebra 
maps $\HH(\iota)\colon\, \HH(N_\G) \to \HH(G_\G)$ 
and $\gr(\iota)\colon\, \gr(N_\G) \to \gr(G_\G)$, which 
commute with the natural surjections from the 
holonomy to the associated graded Lie algebras.   
Passing to derived Lie subalgebras, and tensoring 
with $\Q$, we obtain the following commuting diagram:
\begin{equation}
\label{eq:lie cd}
\xymatrixcolsep{42pt}
\xymatrix{
\HH' (N_\G)\otimes \Q \ar[r]^{\HH'(\iota)\otimes \Q} 
\ar[d]^{\Psi_{N_{\G}} \otimes \Q} 
& \HH' (G_\G)\otimes \Q \ar[d]^{\Psi_{G_{\G}} \otimes \Q} \\
\gr' (N_\G)\otimes \Q  \ar[r]^{\gr'(\iota)\otimes \Q} 
& \gr' (G_\G)\otimes \Q
}
\end{equation}
The vertical arrows are isomorphisms, by the $1$-formality 
of $G_\G$ and $N_\G$, insured by \cite{KM} and 
Proposition~\ref{prop:bb formal}, respectively.  
The bottom arrow is an isomorphism, by 
Theorem \ref{thm:gr chen bb}\ref{bb1}.  
Hence, the top arrow, $\HH'(\iota)\otimes \Q=\HH'_\Q(\iota)$, 
is also an isomorphism.  
\end{proof}

\section{Cohomology ring}
\label{sect:cohoring}

In this section, we exploit the $1$-formality property of a 
finitely presented Bestvina-Brady group $N_\G$, in order 
to give a purely combinatorial description of  $H^{\le 2}(N_\G,\Q)$. 

\subsection{Homology of $N_\G$}
\label{subsec:coho groups}

Fix a coefficient field $\k$.  The extension 
$1\to N_\G\stackrel{\iota}{\to} G_\G \stackrel{\nu}{\to} \Z\to 0$ 
 defines a natural $\k\Z$-module 
structure on $H_*(N_\G,\k)$.  The next result gives a 
combinatorial description of this structure.  
(See also \cite{JM}, \cite{LeS} for related computations.)

\begin{prop}
\label{prop=hbb}
Let $\G$ be a finite graph, with flag complex $\Delta_\G$.  
For each $r>0$, we have 
an isomorphism of $\k\Z$-modules, 
\[
H_r(N_\G,\k) \cong (\k\Z)^{\dim \widetilde{H}_{r-1}(\Delta_\G, \k)} 
\oplus ({}_\epsilon \k)^{\dim B_{r-1}(\Delta_\G, \k)},
\]
where ${}_\epsilon \k$ denotes the trivial $\k\Z$-module $\k$, 
and $B_{\bullet}(\Delta_\G, \k)$ are the simplicial boundaries.
\end{prop}

\begin{proof}
By Shapiro's Lemma, $H_*(N_{\G},\k)=H_*(\k\Z\otimes _{\Z G_{\G}} 
\widetilde{C}_{\bullet} )$, where $\widetilde{C}_{\bullet}$ 
is the equivariant chain complex from Proposition 
\ref{prop:chain complex}, and the change of rings 
$\Z G_\G \to \k \Z$ is induced by $\nu \colon\, G_\G \to \Z$.
Write $P:=\k\Z=\k[\tau^{\pm 1}]$.  Using (\ref{eq:bdry map}), 
we find that 
\begin{equation}
\label{eq=hbb1}
H_r(N_{\G},\k)=\frac{P\otimes Z_{r-1}(\Delta_\G,\k)}{(\tau-1)P
\otimes B_{r-1}(\Delta_\G,\k)},
\end{equation}
as modules over $P$, where  $Z_{\bullet}(\Delta_\G,\k)$ denotes the 
reduced simplicial cycles.  

Set $B:=B_{r-1}(\Delta_\G,\k)$, $Z:=Z_{r-1}(\Delta_\G,\k)$, 
and $H:=\widetilde{H}_{r-1}(\Delta_\G,\k)$. 
It is straightforward to check that the natural maps 
$\k \otimes B \inj P\otimes Z$ and  $P\otimes Z\surj P\otimes H$ 
give rise to the following split exact sequence of $P$-modules:
\begin{equation}
\label{eq=hbb2}
\xymatrix{
0\ar[r]& {}_\epsilon \k \otimes B \ar[r]& 
\DS{\frac{P\otimes Z}{(\tau-1)P \otimes B}}
\ar[r]& P\otimes H \ar[r]& 0
}.
\end{equation}
The conclusion follows by putting together (\ref{eq=hbb1}) 
and  (\ref{eq=hbb2}).
\end{proof}

\subsection{Cohomology ring in low degrees}
\label{subsec=cohring}

Recall from Section \ref{sec:salvetti} that 
$K_{\G}=K(G_\G,1)$ is a subcomplex of the standard 
torus $(S^1)^{\V}$.  This readily implies that the cohomology 
ring $H^*(G_\G,\k)$ is the quotient of the exterior $\k$-algebra 
on $\V$ by the ideal generated by the monomials $vw$ corresponding 
to non-edges $\{v,w\}\in \overline{\E}$. 

\begin{lemma}
\label{cor=1gen}
If $\pi_1(\Delta_{\G})=0$, then the following hold. 
\begin{enumerate}
\item \label{s1}
The cup-product map 
$\cup_{N_\G}\colon\, \bigwedge^2 H^1(N_\G,\Q)\to H^2(N_\G,\Q)$ 
is surjective.
\item \label{s2}
The map $\iota^* \colon\,   H^2(G_\G,\Q)\to  H^2(N_\G,\Q)$ 
is surjective.
\end{enumerate}
\end{lemma}

\begin{proof}
\ref{s1} It is enough to show that 
$\dim_\Q H_2(N_\G,\Q)=\dim_\Q \im (\nabla_{N_\G})$. 
We know from Proposition \ref{prop=hbb} that 
$\dim H_2(N_\G,\Q)= \dim Z_1(\Delta_\G,\Q)=\abs{\E}-\abs{\V}+1$, 
since $\Delta_\G$ is simply-connected.  On the other hand,
\begin{eqnarray*}
\dim \im (\nabla_{N_\G}) & = \dim \bigwedge\nolimits^2 H_1(N_\G,\Q) - 
\dim \HH_2(N_\G)\otimes \Q \\
& = {\abs{\V}-1 \choose 2} - \Big(  {\abs{\V} \choose 2} -\abs{\E}\Big) \\
&= \abs{\E}-\abs{\V}+1, 
\end{eqnarray*}
by the definition of the holonomy Lie algebra,
 Lemma \ref{lem:hhprime}, 
and the injectivity of $\nabla_{G_\G}$.

\ref{s2} Follows from Part \ref{s1}, since we know 
from Lemma \ref{lem:inj map} that $\iota^*\colon\, 
H^1(G_\G,\Q)\to  H^1(N_\G,\Q)$ is surjective.
\end{proof}

\subsection{Proof of Theorem \ref{thm:cohobbintro}}
\label{subsec:proof coho}

Clearly, $\iota^*$ factors through the quotient by the 
ideal generated by $\nu$, since $\nu\iota=0$.  
The isomorphism claim in degree $1$ 
follows immediately from Proposition \ref{lem=alexng}\ref{ng1}.
The surjectivity property in degree $2$ 
is a direct consequence of Lemma \ref{cor=1gen}\ref{s2}.
We are left with proving that 
\[
\ker\left( \iota^*\colon\, H^2(G_\G,\Q)\to H^2(N_\G,\Q)\right) \subset 
\im \left(\cdot \nu\colon\, H^1(G_\G,\Q)\to H^2(G_\G,\Q)\right) .
\]
By dualizing, it is enough to check that the inclusion 
\[
\ker\big( \mu_{\nu}^{\top} \circ \nabla_{G_\G} \colon\, 
H_2(G_\G,\Q)\to H_1(G_\G,\Q)\big) \subset 
\im \big(\iota_* \colon\, H_2(N_\G,\Q)\to H_2(G_\G,\Q)\big) ,
\]
holds, where  $\mu_{\nu}^{\top}$ is the transpose of 
$\mu_\nu\colon\, H^1(G_\G,\Q) \to  \bigwedge^2 
H^1(G_\G,\Q)$, the right-multipli\-cation by $\nu$. 
It follows from Proposition \ref{lem=alexng}\ref{ng1} 
that $\ker \big( \mu_{\nu}^{\top}\big) = 
\im \big( \bigwedge^2 \iota_* \big)$. 
Hence, $\ker\big( \mu_{\nu}^{\top} \circ \nabla_{G_\G}\big)=
 \nabla_{G_\G}^{-1} \big( \im \big( \bigwedge^2 \iota_* \big) 
 \big)$.  Now recall from Lemma \ref{lem:hhprime} that 
 $\bigwedge^2 \iota_*$ induces an isomorphism 
\begin{equation}
\label{eq=rl2}
\bigwedge\nolimits^2 \iota_* \colon\, 
\bigwedge\nolimits^2 H_1(N_\G,\Q)/
\im (\nabla_{N_\G})  \stackrel{\:\cong\:}{\longrightarrow}
 \bigwedge\nolimits^2 H_1(G_\G,\Q)/
\im (\nabla_{G_\G}) .
\end{equation}
The desired inclusion follows at once from 
(\ref{eq=rl2}) and diagram (\ref{eq:nabla}).

This finishes the proof of Theorem \ref{thm:cohobbintro} from 
the Introduction.  In \cite{LeS},  Leary and Saadeto\u{g}lu, 
using a different approach, obtain a similar description of 
the truncated cohomology ring $H^{\le r} (N_\Gamma)$, 
in the situation when $\widetilde{H}_{< r} (\Delta_{\Gamma})=0$.

\section{Characteristic and resonance varieties}
\label{sect:res}

In previous work \cite{PS-artin}, \cite{DPS}, we determined the 
resonance and characteristic varieties of right-angled Artin groups. 
In this section, we do the same for finitely presented 
Bestvina-Brady groups, thus proving 
Theorem \ref{thm:resbbintro} from the Introduction.

\subsection{Jumping loci for $G_\G$}
\label{subsec:res artin}

Let $\G=(\V,\E)$ be a finite graph, 
and let $G_\G$ be the corresponding right-angled Artin group.  
Write $H_{\V}=H^1(G_\G, \C)$ and $\TT_\V=\Hom(G_\G,\C^*)$.  
If $\W$ is a subset of $\V$, write $H_{\W}$ and $\TT_\W$ 
for the corresponding coordinate subspaces, respectively, subtori. 

\begin{theorem}[(\cite{PS-artin}, \cite{DPS})]
\label{thm:res artin}
Let $\G$ be a finite graph.  Then:
\[ 
\label{eq:res artin}
\RR(G_{\G}) = \bigcup_{\stackrel{\W\subset \V}{
\G_{\W} \textup{\scriptsize{ disconnected}}}} H_{\W}
\qquad\mbox{and}\qquad 
\VV(G_{\G}) = \bigcup_{\stackrel{\W\subset \V}{
\G_{\W} \textup{\scriptsize{ disconnected}}}} \TT_{\W}.
\] 
\end{theorem}

It follows that the irreducible components of the above 
varieties are indexed by the subsets $\W\subset \V$, 
maximal among those for which $\G_{\W}$ is disconnected. 
 In particular, if $\G$ is disconnected, 
then $\RR(G_{\G})=H_{\V}$ and $\VV(G_{\G})=\TT_{\V}$. 

\begin{example}
\label{ex:artin tree res}
Let $\G$ be a tree on $n>2$ vertices.  Then  
$\RR(G_\G)$ is a union of coordinate hyperplanes, 
one for each cut point (that is, non-extremal vertex) 
of $\G$.  In particular, the irreducible components 
of $\RR(G_\G)$ are in general position.  
\end{example}

\subsection{A map between jumping loci}
\label{subsec:bb res} 
Let $N_{\G}=\ker(\nu\colon\, G_\G \surj \Z)$, and let 
$\iota\colon\, N_\G\inj G_\G$ be the natural inclusion. 

Assuming $\G$ is connected, 
we infer from Proposition \ref{lem=alexng}\ref{ng1} 
that $\iota$ induces a vector space epimorphism 
$\iota^*\colon\, H^1(G_\G,\C) \surj H^1(N_\G,\C)$. 
Identifying $H^1(G_\G,\C)$ with $\C^{\V}$, the  kernel of 
$\iota^*$ gets identified with the diagonal line.  

Similarly,  $\iota$ induces  an algebraic group  
epimorphism $\iota^*\colon\, \TT_{G_\G} \surj \TT_{N_\G}$.  
Identifying $\TT_{G_\G}$ with $(\C^*)^{\V}$, the  kernel of 
$\iota^*$ gets identified with the diagonal $1$-dimensional 
subtorus. 

\begin{lemma}
\label{lem:res bb artin}
If the flag complex $\Delta_{\G}$ is simply-connected, 
and $\abs{\V}>1$, then 
\begin{enumerate}
\item \label{resbb1}
The map $\iota^*\colon\, H^1(G_\G, \C) \surj H^1(N_\G, \C)$ 
restricts to a surjection  
$\iota^* \colon\, \RR(G_\G)\surj \RR(N_\G)$.  

\item \label{resbb2}
The map $\iota^*\colon\, \TT_{G_\G} \surj \TT_{N_\G}$ 
restricts to a surjection  
$\iota^* \colon\, \VV(G_\G)\surj \VV(N_\G)$.  
\end{enumerate}
\end{lemma}

\begin{proof}
Part \ref{resbb1}. There is an intimate connection 
between the resonance variety and the infinitesimal 
Alexander invariant of a finitely presented group 
$G$.  More precisely, $\RR(G)$ coincides, away from the
origin, with the zero set of the annihilator of
the $\Sym(H_1(G))\otimes\C$-module $\B(G)\otimes \C$. 
This follows from \cite[Lemma 4.2]{DPS}, together with 
the description of radicals of Fitting ideals in terms of 
annihilators, see \cite[pp.~511--513]{E}.

It follows from Lemma \ref{lem:hhprime} that 
$\B(\iota) \otimes \C\colon\, \B (N_\G) \otimes \C\to 
\B (G_\G) \otimes \C$ is an isomorphism of modules 
over $\Sym(H_1(N_\G)) \otimes \C$. Hence, 
\begin{equation}
\label{eq:ann}
\ann ( \B (N_\G) \otimes \C) =  
(\iota_* \otimes \C)^{-1} (\ann ( \B (G_\G) \otimes \C)).
\end{equation}
Taking the complex varieties defined by the ideals on 
both sides of (\ref{eq:ann}) finishes the proof.

Part \ref{resbb2}. 
Similarly, we know from \cite[Lemma 4.5]{DPS} that 
$\VV(G)$ coincides, away from the
origin, with the zero set of the annihilator of
the $\C H_1(G)$-module $B(G)\otimes \C$. 
By Proposition \ref{lem=alexng}\ref{ng2}, the map 
$B(\iota) \otimes \C\colon\, B(N_{\G}) \otimes \C\to 
B(G_{\G}) \otimes \C$ is an isomorphism of modules 
over $\C H_1(N_{\G})$. As above, we conclude by taking 
the complex varieties defined by the annihilators of 
these two modules. 
\end{proof}

\begin{lemma}
\label{lem:inj}
Suppose $\W$ is a proper subset of the vertex set $\V$ of $\G$. 
Then the restriction of $\iota^* \colon\, H^1(G_\G, \C)\to H^1(N_\G, \C)$  
to the coordinate subspace $H_{\W}$ is injective.  Similarly, 
the restriction of $\iota^*\colon\, \TT_{G_\G} \to \TT_{N_\G}$  
to the coordinate subtorus $\TT_{\W}$ is injective. 
\end{lemma}

\begin{proof}
Let $\{e_v\}_{v \in \V}$ be the standard basis for the vector 
space $H^1(G_{\G},\C)=\C^{\V}$.  Suppose 
$\iota^*\big(\sum_{w\in \W} c_w  e_w\big)=0$.  
Then $ \sum_{w\in \W} c_w  e_w =  c \sum_{v\in \V} e_v $, 
for some scalar $c \in \C$.  Picking $v\in \V\setminus \W$, 
and comparing the coefficient of $e_v$ on both sides of 
this equation, we see that $c=0$. This finishes the proof of 
the first claim; the proof of the second claim follows 
along the same lines.
\end{proof}

For a subset $\W\subset \V$, let $H'_{\W}$ denote the 
subspace $\iota^*(H_{\W})\subset H^1(N_\G, \C)$, and let 
$\TT'_{\W}$ denote the 
subtorus $\iota^*(\TT_{\W})\subset \TT_{N_\G}$.

\begin{lemma}
\label{lem:intersection}
Suppose  $\W_1$ and $\W_2$ are two 
subsets of $\V$, of size at most $\abs{\V}-2$.  If  
$H'_{\W_1}\subset H'_{\W_2}$, or 
$\TT'_{\W_1}\subset \TT'_{\W_2}$, 
 then $\W_1 \subset \W_2$.  
\end{lemma}

\begin{proof}
Assume there is a vertex $v_1\in \W_1\setminus \W_2$.  
Since $\abs{\W_2}\le \abs{\V}-2$, there must be another 
vertex $v_2\in\V\setminus  \W_2$, distinct from $v_1$.  
Suppose $H'_{\W_1}\subset H'_{\W_2}$.  Then $e_{v_1} = 
\sum_{v\in \W_2} c_v  e_v + c \sum_{v\in \V} e_v$.  
Comparing coefficients of $e_{v_2}$ on both sides, 
we find $c=0$; hence $e_{v_1}\in H_{\W_2}$, 
a contradiction.  The case $\TT'_{\W_1}\subset \TT'_{\W_2}$ 
is treated similarly.
\end{proof}

\subsection{Proof of Theorem \ref{thm:resbbintro}}
\label{subsec:proof jump}

It follows from Theorem \ref{thm:res artin} and 
Lemma \ref{lem:res bb artin} that 
\begin{equation} 
\label{eq:jumpbb}
\RR(N_{\G}) = \bigcup_{\W} H'_{\W}
\qquad\mbox{and}\qquad 
\VV(N_{\G}) = \bigcup_{\W} \TT'_{\W}, 
\end{equation} 
where, in both cases, the union is taken over 
all subsets $\W\subset \V$, maximal 
among those for which $\Gamma_{\W}$ is disconnected. 
Lemma~\ref{lem:inj} guarantees that 
$H'_{\W}\subset H^1(N_\G,\C)$ is a 
vector subspace of dimension $\abs{\W}$, and 
$\TT'_{\W}\subset \TT_{N_\G}$ is a subtorus of 
dimension $\abs{\W}$. 

First assume $\kappa(\G)=1$, that is, $\G$ has a cut point.  
This means there is a $v\in \V$ such that 
$\G_{\V\setminus v}$ is disconnected.  A dimension 
count shows that  $\RR(N_\G)=H^1(N_{\G},\C)$ and 
$\VV(N_\G)=\TT_{N_{\G}}$. 

Now assume $\kappa(\G)>1$.  
We infer from  Lemma~\ref{lem:intersection} that 
(\ref{eq:jumpbb}) gives indeed 
the irreducible decompositions of the respective 
varieties.  This ends the proof of 
Theorem \ref{thm:resbbintro}.

\begin{example}
\label{ex:bb tree res}
Let $\G$ be a tree on $n>2$ vertices.  Then  $\kappa(\G)=1$, 
and so $\RR(N_\G)=H^1(N_\G,\C)=\C^{n-1}$, and 
$\VV(N_\G)=\TT_{N_\G}=(\C^*)^{n-1}$ (this computation 
also follows from the fact that $N_\G=F_{n-1}$).
\end{example}

\section{Comparison with Artin groups and arrangement groups}
\label{sect:artin-arr}

In this section, we use the methods developed above 
to compare the Bestvina-Brady groups to two other classes 
of groups: Artin groups and arrangement groups.

\subsection{Extra-special triangulations}
\label{subsec:extraspecial}

Recall we defined a triangulation $\Delta$ of the disk $D^2$ 
to be  {\em special}\/ if $\Delta$ can be 
obtained from a triangle by adding one triangle at 
a time, along a unique boundary edge.  

 \begin{lemma}
 \label{lem:special}
 Let $\Delta$ be a special triangulation of the $2$-disk, 
 and $\G=(\V,\E)$ its $1$-skeleton. 
 Then $\RR(N_\G)$ is a proper subset of $H^1(N_\G,\C)$.  
 \end{lemma}
 
 \begin{proof}
Recall from Lemma \ref{lem:special triang}\ref{sp2} 
that $\Delta_\G=\Delta$; in particular,  Theorem 
\ref{thm:resbbintro} applies.  It is also  readily seen that 
 the graph $\G$ has no cut points, i.e, $\kappa(\G)>1$. 
 Thus, $\RR(N_\G)\subsetneq H^1(N_\G,\C)$.  
 \end{proof}

\begin{definition}
\label{def:triang}
 A triangulation of $D^2$ is called {\em extra-special}\/ 
 if it is obtained from a special triangulation, 
by adding one triangle along each boundary edge. 
(See Figure \ref{fig:two triangulations}.)
\end{definition}

\begin{figure}
\subfigure{%
\begin{minipage}[t]{0.25\textwidth}
\setlength{\unitlength}{0.5cm}
\begin{picture}(4,5)(-0.5,-0.5)
\put(0,0){\line(1,0){4}}
\put(0,4){\line(1,0){4}}
\put(0,0){\line(0,1){4}}
\put(4,0){\line(0,1){4}}
\put(0,0){\line(1,1){4}}
\put(6,2){$\leadsto$}
\end{picture}
\end{minipage}
}
\subfigure{%
\begin{minipage}[t]{0.25\textwidth}
\setlength{\unitlength}{0.3cm}
\begin{picture}(4,5.2)(-5.5,-2.2)
\put(0,0){\line(1,0){4}}
\put(0,4){\line(1,0){4}}
\put(0,0){\line(0,1){4}}
\put(4,0){\line(0,1){4}}
\put(0,0){\line(1,1){4}}
\put(-2,2){\line(1,-1){4}}
\put(-2,2){\line(1,1){4}}
\put(6,2){\line(-1,-1){4}}
\put(6,2){\line(-1,1){4}}
\end{picture}
\end{minipage}
}
\caption{\textsf{Building an extra-special 
triangulation of the disk}}
\label{fig:two triangulations}
\end{figure}

If $\Delta$ is extra-special, more 
can be said about the resonance variety of $N_\G$.  
By definition, $\Delta$ is obtained by 
attaching triangles to the boundary edges  of a special 
triangulation of the disk. Denote by   
$(e_1,\dots ,e_r)$ the circuit formed by these edges, 
and write $\W_i=\V\setminus e_i$. 
Note that each edge $e_i$ is a minimal cut set of $\G$; 
hence, $H'_{\W_i}$ is an irreducible component 
of $\RR(N_{\G})$, for each $i=1,\dots, r$. 

\begin{lemma}
\label{lem:extraspecial}
Let $\G$ be the $1$-skeleton of an extra-special triangulation 
$\Delta$ of $D^2$. Then the subspace $\bigcap_{i=1}^{r} 
H'_{\W_i}$ has codimension $r-1$ in $H^1(N_\G,\C)$. 
In particular,  $\bigcap_{i=1}^{r} H'_{\W_i} \ne 0$. 
\end{lemma}

\begin{proof} 
We claim that
\begin{equation}
\label{eq:extra}
\iota^* \Big(\bigcap_{i=1}^r  H_{\W_i}\Big)=\bigcap_{i=1}^r  H'_{\W_i}.
\end{equation}

The inclusion $\subseteq$ is clear.  The reverse inclusion 
is proved by induction on $s$ ($0<s<r$), with the case $s=1$ 
being obvious. Set 
$P_{k}=\bigcap _{i=1}^{k} H_{\W_i}$, $Q_{k}= H_{\W_k}$, 
$P'_{k}=\bigcap _{i=1}^{k} H'_{\W_i}$, and $Q'_{k}= H'_{\W_k}$. 
We then have the commuting diagram
\begin{equation}
\label{eq:cdp}
\xymatrix{
0\ar[r] & 
P_{s}\cap Q_{s+1} \ar^{\iota^*}[d] \ar[r] 
& P_{s}\oplus Q_{s+1} \ar^{\iota^*}[d] \ar[r] 
& P_{s}+ Q_{s+1}  \ar^{\iota^*}[d] \ar[r] & 0  \\
0\ar[r] 
&P'_{s}\cap Q'_{s+1} \ar[r] 
&P'_{s}\oplus Q'_{s+1} \ar[r] 
&P'_{s} + Q'_{s+1}  \ar[r] 
& 0}
\end{equation}
The middle arrow is an isomorphism, by Lemma \ref{lem:inj} 
and the induction hypothesis.  Clearly, the right arrow 
is an epimorphism.  
Note that $P_{s}+ Q_{s+1}$ is a subspace of 
\[
H_{\W_1\cap \cdots \cap \W_s} + H_{\W_{s+1}}= 
H_{(\W_1\cap \cdots \cap \W_s)\cup \W_{s+1}}.
\]
Since $(e_1\cup\cdots \cup e_{s})\cap e_{s+1}\ne \emptyset$,  
this is a proper subspace of $H_{\V}$.  
Thus, the right arrow in diagram (\ref{eq:cdp}) is injective, 
again by Lemma \ref{lem:inj}.  Applying the $5$-Lemma 
finishes the proof of the claim.

From (\ref{eq:extra}), we see that 
\[
\label{eq:nontransverse}
\codim \bigcap_{i=1}^r  H'_{\W_i}= 
( \abs{\V} -1 ) - \dim \bigcap_{i=1}^r  H_{\W_i}=r-1. 
\]
Finally, if $ \bigcap_{i=1}^r  H'_{\W_i}=0$, then $r=\abs{\V}$. 
But clearly $\abs{\V}\ge 2r$. 
\end{proof}

\subsection{Artin groups}
\label{subsec:artin}
A {\em weighted graph}\/ is a graph $\G=(\V,\E)$ 
endowed with a function $m\colon\, \E\to \Z$ that assigns 
to each edge $e$ an integer $m(e)\ge 2$.  
Such a weighted graph $(\G,m)$  determines 
an Artin group (of finite type),
\[
\label{eq:artin}
G_{\G,m}= \langle v \in \V \mid \pi_m(v,w) =\pi_m(w,v) 
\mbox{  if } \{v,w\} \in \E\rangle,
\]
where $\pi_m(v,w)=vwv\cdots$ has length $m(\{v,w\})$.  
When all edge weights are equal to $2$, this is simply the 
right-angled Artin group $G_{\G}$. 

Associated to a weighted graph as above there is an 
ordinary (unlabeled) graph $\widetilde{\G}=
({\widetilde  \V},{\widetilde  \E})$, called 
the ``odd contraction" of $(\G,m)$, see  \cite[\S10.9]{DPS}.  
First define $\Gamma_{\rm{odd}}$ to be the unlabeled 
graph with vertex set $\V$ and edge set 
$\{e \in \E \mid  m (e)~  \mbox{is odd}\}$.
Then define ${\widetilde  \V}$ to be the set of connected 
components of $\Gamma_{\rm{odd}}$, with two distinct 
components determining an edge $\{c,c'\} \in  {\widetilde \E}$ 
if and only if there exist vertices $v\in c$ 
and $v'\in c'$ which are connected by an edge in $\E$.

\begin{prop}
\label{prop:special bb artin}
Let $\G$ be the $1$-skeleton of an extra-special triangulation 
of $D^2$.  Then the Bestvina-Brady group $N_\G$ is not 
isomorphic to any Artin group.
\end{prop}

\begin{proof}
Suppose $N_\G$ is isomorphic to an Artin group 
$G_{\G',m}$.  Let $\widetilde{\G}'$ be the odd 
contraction of $\G'$.   Lemma~10.11 from \cite{DPS} guarantees 
that the respective Malcev Lie algebras, $M_{G_{\G',m}}$  
and $M_{G_{\widetilde{\G}'}}$, are filtered isomorphic. 
On the other hand, we know from Theorem~16.10 from \cite{KM} 
that both Artin groups, $G_{\G',m}$ and $G_{\widetilde{\G}'}$, 
are $1$-formal. Passing to associated graded Lie algebras, 
we obtain that 
$\HH_{\Q} (G_{\G',m}) \cong \HH_{\Q} (G_{\widetilde{\G}'})$, 
as graded Lie algebras.  Hence, 
$\HH_{\Q} (N_{\G}) \cong \HH_{\Q} (G_{\widetilde{\G}'})$, 
as graded Lie algebras. This implies the existence of 
an ambient isomorphism
\[
\RR(N_{\G}) \cong \RR(G_{\widetilde{\G}'}).
\]
From Lemma \ref{lem:special}, 
we know that $\RR(N_\G)\subsetneq H^1(N_\G,\C)$.  
By Theorem \ref{thm:res artin}, 
this forces $\widetilde{\G'}$ to be connected. 

Let $\G=(\V,\E)$, and write $v=\abs{\V}$, $e=\abs{\E}$. 
Similarly, let $\widetilde{\G}'=(\V',\E')$, and write $v'=\abs{\V'}$, 
$e'=\abs{\E'}$. We claim $v'=e'+1$, and thus, $\widetilde{\G}'$ 
is a tree.

Note that $v'=b_1(G_{\widetilde{\G'}})=b_1(N_\G)=v-1\ge 5$. 
Moreover, ${v' \choose 2}-e'=\rank \HH_2(G_{\widetilde{\G}'})= 
\rank \HH_2(N_\G)$.  We also know that $\rank \HH_2(N_\G)= 
\rank \HH_2(G_\G)$, by Lemma \ref{lem:hhprime}.  Since 
$\rank \HH_2(G_\G)= {v \choose 2}-e$, we conclude that 
$e'=e-v+1$. From Lemma \ref{lem:special triang}\ref{sp1}, 
we know $2v-e=3$; hence, $v'=e'+1$, as claimed. 

By the discussion from Example \ref{ex:artin tree res}, the 
components of $\RR(N_\G)=\RR(G_{\widetilde{\G}'})$ 
must intersect transversely.  This contradicts 
Lemma \ref{lem:extraspecial}.  
\end{proof}

\subsection{Arrangement groups}
\label{subsec:arr}

Another widely studied class of groups are the fundamental 
groups of complements of complex hyperplane arrangements; 
see for instance \cite{Su} and references therein.  
Bestvina-Brady groups associated to simply-connected 
flag complexes share some common features with 
arrangement groups.  Indeed, if $G$ is a group in 
either class, then:
\begin{itemize}
\item \label{1}
$G$ admits a finite presentation, with 
commutator relators only;
\item \label{2}
$G$ is $1$-formal;
\item  \label{3}
$\RR(G)$ is a union of linear subspaces.
\end{itemize}

There is a rather striking similarity between 
Bestvina-Brady groups associated to complete 
multipartite graphs and the fundamental groups 
of complements of ``decomposable" arrangements.  
Indeed, if $G$ is a group in either class, then the 
derived Lie algebra of $\gr(G)$ is isomorphic to the 
derived Lie algebra of a product of free groups, and 
similarly for the derived Chen Lie algebra.  For 
Bestvina-Brady groups, this was noted in 
Remark \ref{rem:bb}, while for 
decomposable arrangement groups, this is 
proved (by completely different methods) in 
Theorems 2.4 and 6.2 from \cite{PS-decomp}. 

Even so, there are many finitely presented 
Bestvina-Brady groups which are not arrangement groups. 

\begin{prop}
\label{prop:special bb arr}
Let $\G$ be the $1$-skeleton of an extra-special triangulation 
of $D^2$.  Then the Bestvina-Brady group $N_\G$ is not 
isomorphic to any arrangement group.  
\end{prop}

\begin{proof}
If $G$ is an arrangement group, then any two  
components of $\RR(G)$ intersect only at $0$, 
see \cite{LY}.  By Lemma \ref{lem:extraspecial}, 
this cannot happen for $N_\G$. 
\end{proof}

Propositions \ref{prop:special bb artin} and 
\ref{prop:special bb arr} together yield 
Theorem \ref{thm:artinbbintro} from the Introduction.

\begin{ack}
This work was started while the two authors were 
attending the program ``Hyperplane Arrangements 
and Applications'' at the Mathematical Sciences 
Research Institute in Berkeley, California, in Fall, 2004.  
A preliminary version was presented by the second 
author in a talk delivered at the MSRI workshop on 
``Combinatorial Aspects of Hyperplane Arrangements''
(November 1--5, 2004). We thank MSRI for providing an 
excellent research environment.  The paper was completed 
while the first author was visiting Northeastern University, 
in Spring, 2006. He thanks the Northeastern Mathematics 
Department for its support and hospitality during his stay.
\end{ack}

\affiliationone{%
Stefan Papadima\\
Inst.~of Math.~Simion Stoilow \\
P.O. Box 1-764\\
RO-014700 Bucharest\\ 
Romania \\
\email{Stefan.Papadima@imar.ro}}
\affiliationtwo{%
Alexander~I.~Suciu \\
Department of Mathematics\\
Northeastern University\\
Boston, MA 02115\\ 
USA\\
\email{a.suciu@neu.edu}}

\end{document}